\documentclass{amsart}

\usepackage{cite}

\numberwithin{equation}{section}

\newtheorem{lemma}{Lemma}[section]
\newtheorem{definition}{Definition}[section]
\newtheorem{theorem}{Theorem}[section]

\begin{document}

\title[The well-poised Bailey lemma]
{Extensions of the well-poised and elliptic well-poised Bailey lemma}

\author{S.~Ole Warnaar}\thanks{Work supported by the Australian
Research Council}
\address{Department of Mathematics and Statistics,
The University of Melbourne, Vic 3010, Australia}
\email{warnaar@ms.unimelb.edu.au}



\begin{abstract}
We establish a number of extensions of the well-poised Bailey lemma
and elliptic well-poised Bailey lemma.
As application we prove some new transformation formulae
for basic and elliptic hypergeometric series, and
embed some recent identities of Andrews, Berkovich and Spiridonov
in a well-poised Bailey tree.
\end{abstract}

\dedicatory{Dedicated to Tom Koornwinder on the occasion of his sixtieth 
birthday}

\maketitle

\section{Introduction}
In a recent paper \cite{Andrews01} Andrews introduced a Bailey-type
lemma for well-poised (WP) basic series.
Together with Berkovich this led him to discover many new 
transformation formulae for basic hypergeometric series \cite{AB02}.
Shortly after \cite{Andrews01} and \cite{AB02} appeared, 
Spiridonov \cite{Spiridonov02} pointed out that part of the programme 
carried out by Andrews and Berkovich can be formulated at the level
of elliptic functions, leading to new results for elliptic or modular 
hypergeometric series.

In the present paper we show that the work of Andrews, Berkovich and
Spiridonov admits many further extensions. In particular we
will show that Andrews' binary WP Bailey tree can be enhanced to
yield a tree with six-fold branching, and that Spiridonov's elliptic 
WP Bailey chain can be upgraded to a trinary tree.
As a consequence, all of the basic WP Bailey pairs of
\cite{Andrews01,AB02} and all of the elliptic WP Bailey pairs of
\cite{Spiridonov02} become nodes on a single multi-dimensional
basic or elliptic WP Bailey tree. 
In addition, many new WP Bailey pairs result, leading to
new identities for basic and elliptic hypergeometric series.

In Section~\ref{secbasic} we present our extensions to the
basic WP Bailey lemma, Section~\ref{secell} deals with the
elliptic WP Bailey lemma, and Section~\ref{secapl} contains a number
of applications of our result to basic and elliptic hypergeometric series.

\section{The basic well-poised Bailey lemma}\label{secbasic}
Throughout this paper we adopt the standard notation and terminology
for basic hypergeometric series of Gasper and Rahman's book \cite{GR90}.
In particular we use
\begin{equation*}
{_{r+1}\phi_r}\biggl[\genfrac{}{}{0pt}{}
{a_1,\dots,a_{r+1}}{b_1,\dots,b_r};q,z\biggr]
=\sum_{k=0}^{\infty} 
\frac{(a_1,\dots,a_{r+1};q)_k}{(q,b_1,\dots,b_r;q)_k}\, z^k,
\end{equation*}
where 
$(a;q)_n=\prod_{j=0}^{n-1}(1-aq^j)$ is a $q$-shifted factorial and
\begin{equation*}
(a_1,\dots,a_k;q)_n=(a_1;q)_n\dots(a_k;q)_n.
\end{equation*}
Since we will only be dealing with terminating series we do not impose 
the usual condition $|q|<1$.
In view of the frequent use of the term `well-poised' we recall that
a $_{r+1}\phi_r$ series is well-poised if the pairwise
product of numerator and denominator parameters is constant;
$qa_1=a_2b_1=\cdots=a_{r+1}b_r$. If in addition 
$a_2=-a_3=a_1^{1/2}q$, the series is very-well-poised.
We abbreviate such very-well-poised series by
$_{r+1}W_r(a_1;a_4,\dots,a_{r+1};q,z)$. Finally we note our convention that
$n$ is always a nonnegative integer.

\vspace{4mm}

Before discussing Andrews' WP Bailey lemma let us give the
well-poised version of the classical Bailey transform.
\begin{lemma}[WP Bailey transform]\label{lemT}
For $a$ and $k$ indeterminates the following two equations are equivalent:
\begin{subequations}\label{WPBT}
\begin{align}\label{ab}
\beta_n(a,k;q)&=\sum_{r=0}^n \frac{(k/a;q)_{n-r}}{(q;q)_{n-r}}
\frac{(k;q)_{n+r}}{(aq;q)_{n+r}}\,\alpha_r(a,k;q), \\
\alpha_n(a,k;q)&=
\frac{1-aq^{2n}}{1-a}\sum_{r=0}^n\frac{1-kq^{2r}}{1-k}
\frac{(a/k;q)_{n-r}}{(q;q)_{n-r}}
\label{ba} \\
& \qquad \qquad \qquad \qquad \times \frac{(a;q)_{n+r}}{(kq;q)_{n+r}}
\Bigl(\frac{k}{a}\Bigr)^{n-r}\beta_r(a,k;q). \notag
\end{align}
\end{subequations}
\end{lemma}

\begin{proof}
Writing \eqref{WPBT} as
\begin{subequations}
\begin{align}\label{bMa}
\beta_n(a,k;q)&=\sum_{r=0}^n M_{n,r}(a,k;q) \alpha_r(a,k;q), \\
\alpha_n(a,k;q)&=\sum_{r=0}^n \tilde{M}_{n,r}(a,k;q) \beta_r(a,k;q)
\end{align}
\end{subequations}
we get
\begin{align*}
\alpha_n(a,k;q)&=\sum_{r=0}^n \tilde{M}_{n,r}(a,k;q) 
\sum_{s=0}^r M_{r,s}(a,k;q) \alpha_s(a,k;q) \\
&=\sum_{s=0}^n \alpha_s(a,k;q) 
\sum_{r=s}^n \tilde{M}_{n,r}(a,k;q) M_{r,s}(a,k;q).
\end{align*}
We thus need to prove the inverse relation
\begin{equation}\label{MM}
\sum_{s=r}^n \tilde{M}_{n,s}(a,k;q) M_{s,r}(a,k;q)=\delta_{n,r}.
\end{equation}
From the explicit expressions for $\tilde{M}_{n,s}$ and $M_{s,r}$
it follows that
\begin{align*}
\text{LHS}\eqref{MM}
&=\frac{(k;q)_{2r}}{(aq;q)_{2r}}\tilde{M}_{n,r}(a,k;q) 
\,{_6}W_5(kq^{2s};k/a,aq^{n+r},q^{-(n-r)};q,q) \\
&=\frac{(k;q)_{2r}}{(aq;q)_{2r}}\tilde{M}_{n,r}(a,k;q)\,\delta_{n,r}
=\delta_{n,r}.
\end{align*}
Here the ${_6}W_5$ has been summed by a special case of 
Rogers' $q$-Dougall sum \cite[Equation (II.21)]{GR90}.

The companion relation
\begin{equation*}
\sum_{s=r}^n M_{n,s}(a,k;q) \tilde{M}_{s,r}(a,k;q)=\delta_{n,r}
\end{equation*}
can be proved in analogous fashion. In fact, it suffices to note
that given $\alpha_n$ equation \eqref{ba} uniquely fixes 
$\beta_n$ and, conversely, given $\beta_n$ equation \eqref{ab}
uniquely fixes $\alpha_n$.
\end{proof}

After the above preliminaries we review Andrews' work on
the WP Bailey lemma.
Let $\alpha=(\alpha_0,\alpha_1,\dots)$ and
$\beta=(\beta_0,\beta_1,\dots)$.
\begin{definition}
A WP Bailey pair is a pair of sequences
$(\alpha(a,k;q),\beta(a,k;q))$ that satisfies \eqref{WPBT}.
\end{definition}
The adjective `well-poised' derives from the fact that
\begin{equation*}
\frac{(k/a;q)_{n-r}}{(q;q)_{n-r}}\frac{(k;q)_{n+r}}{(aq;q)_{n+r}}
=\frac{(k,k/a;q)_n}{(q,aq;q)_n}
\frac{(kq^n,q^{-n};q)_r}{(aq^{1-n}/k,aq^{n+1};q)_r} 
\Bigl(\frac{aq}{k}\Bigr)^r.
\end{equation*}
When $k$ tends to zero a WP Bailey pair reduces to the
classical Bailey pair as introduced by Bailey \cite{Bailey49}.
For more details we refer to the reviews 
\cite{Andrews84,Andrews01,W01}.

Andrews proved two theorems for constructing WP Bailey pairs.
His first result \cite[Theorem 3]{Andrews01} follows from
Jackson's $_8\phi_7$ sum.
\begin{theorem}\label{thm1}
If $(\alpha(a,k;q),\beta(a,k;q))$ is a WP Bailey pair, then so is the pair
$(\alpha'(a,k;q),\beta'(a,k;q)$ given by
\begin{align*}
\alpha'_n(a,k;q)&=\frac{(b,c;q)_n}{(aq/b,aq/c;q)_n}
\Bigl(\frac{k}{m}\Bigr)^n \alpha_n(a,m;q), \\
\beta'_n(a,k;q)&=\frac{(mq/b,mq/c;q)_n}
{(aq/b,aq/c;q)_n}\sum_{r=0}^n \frac{1-mq^{2r}}{1-m}
\frac{(b,c;q)_r}{(mq/b,mq/c;q)_r} \\
& \qquad \qquad \qquad \qquad \qquad \times
\frac{(k/m;q)_{n-r}}{(q;q)_{n-r}}\frac{(k;q)_{n+r}}{(mq;q)_{n+r}} 
\Bigl(\frac{k}{m}\Bigr)^r \beta_r(a,m;q),
\end{align*}
where $m=bck/aq$.
\end{theorem}
Andrews' second result \cite[Theorem 4]{Andrews01} follows from
the $q$-Pfaff--Saalsch\"utz summation.
\begin{theorem}\label{thm2}
If $(\alpha(a,k;q),\beta(a,k;q))$ is a WP Bailey pair, then so is the pair
$(\alpha'(a,k;q),\beta'(a,k;q)$ given by
\begin{align*}
\alpha'_n(a,k;q)&=\frac{(m;q)_{2n}}{(k;q)_{2n}}
\Bigl(\frac{k}{m}\Bigr)^n \alpha_n(a,m;q), \\
\beta'_n(a,k;q)&=\sum_{r=0}^n \frac{(k/m;q)_{n-r}}{(q;q)_{n-r}}
\Bigl(\frac{k}{m}\Bigr)^r \beta_r(a,m;q),
\end{align*}
where $m=a^2q/k$.
\end{theorem}
By combining the above two theorems it follows that each 
WP Bailey pair gives rise to a binary tree of WP Bailey pairs.
Andrews coined this the WP Bailey tree.

We will now show that the Bailey tree admits many additional branches
thanks to the following four theorems.
\begin{theorem}\label{thm2b}
If $(\alpha(a,k;q),\beta(a,k;q))$ is a WP Bailey pair, then so is the pair
$(\alpha'(a,k;q),\beta'(a,k;q)$ given by
\begin{align*}
\alpha'_n(a,k;q)&=
\frac{1-\sigma k^{1/2}}{1-\sigma k^{1/2}q^n}
\frac{1+\sigma m^{1/2}q^n}{1+\sigma m^{1/2}}
\frac{(m;q)_{2n}}{(k;q)_{2n}}
\Bigl(\frac{k}{m}\Bigr)^n \alpha_n(a,m;q), \\
\beta'_n(a,k;q)&=\frac{1-\sigma k^{1/2}}{1-\sigma k^{1/2}q^n}
\sum_{r=0}^n \frac{1+\sigma m^{1/2}q^r}{1+\sigma m^{1/2}}
\frac{(k/m;q)_{n-r}}{(q;q)_{n-r}}
\Bigl(\frac{k}{m}\Bigr)^r \beta_r(a,m;q),
\end{align*}
where $m=a^2/k$ and $\sigma\in\{-1,1\}$.
\end{theorem}
The freedom in the choice of $\sigma$ simply reflects the fact that
the above expressions are invariant under the simultaneous 
negation of $k^{1/2}$, $m^{1/2}$ and $\sigma$.

\begin{theorem}\label{thm3}
If $(\alpha(a,k;q),\beta(a,k;q))$ is a WP Bailey pair, then so is the pair
$(\alpha'(a,k;q),\beta'(a,k;q)$ given by
\begin{align*}
\alpha'_n(a^2,k;q^2)&=\alpha_n(a,m;q), \\
\beta'_n(a^2,k;q^2)&=\frac{(-mq;q)_{2n}}{(-aq;q)_{2n}}
\sum_{r=0}^n \frac{1-mq^{2r}}{1-m}
\frac{(k/m^2;q^2)_{n-r}}{(q^2;q^2)_{n-r}} \\
&\qquad \qquad \qquad \qquad \quad \times
\frac{(k;q^2)_{n+r}}{(m^2q^2;q^2)_{n+r}}
\Bigl(\frac{m}{a}\Bigr)^{n-r} \beta_r(a,m;q),
\end{align*}
where $m=k/aq$.
\end{theorem}

\begin{theorem}\label{thm4}
If $(\alpha(a,k;q),\beta(a,k;q))$ is a WP Bailey pair, then so is the pair
$(\alpha'(a,k;q),\beta'(a,k;q)$ given by
\begin{align*}
\alpha'_n(a^2,k;q^2)&=q^{-n}\frac{1+aq^{2n}}{1+a}\,\alpha_n(a,m;q), \\
\beta'_n(a^2,k;q^2)&=
q^{-n}\frac{(-mq;q)_{2n}}{(-a;q)_{2n}} 
\sum_{r=0}^n \frac{1-mq^{2r}}{1-m}
\frac{(k/m^2;q^2)_{n-r}}{(q^2;q^2)_{n-r}} \\
&\qquad \qquad \qquad \qquad \qquad \times
\frac{(k;q^2)_{n+r}}{(m^2q^2;q^2)_{n+r}}
\Bigl(\frac{m}{a}\Bigr)^{n-r}\beta_r(a,m;q),
\end{align*}
where $m=k/a$.
\end{theorem}

\begin{theorem}\label{thm5}
If $(\alpha(a,k;q),\beta(a,k;q))$ is a WP Bailey pair, then so is the pair
$(\alpha'(a,k;q),\beta'(a,k;q)$ given by
\begin{align*}
\alpha'_{2n}(a,k;q)&=\alpha_n(a,m;q^2), \qquad  
\alpha'_{2n+1}(a,k;q)=0, \\
\beta'_n(a,k;q)&=\frac{(mq;q^2)_n}{(aq;q^2)_n}
\sum_{r=0}^{\lfloor n/2\rfloor}
\frac{1-mq^{4r}}{1-m}
\frac{(k/m;q)_{n-2r}}{(q;q)_{n-2r}} \\
& \qquad\qquad\qquad\qquad\quad \times
\frac{(k;q)_{n+2r}}{(mq;q)_{n+2r}}
\Bigl(-\frac{k}{a}\Bigr)^{n-2r} 
\beta_r(a,m;q^2),
\end{align*}
where $m=k^2/a$.
\end{theorem}

The Theorems~\ref{thm3} and \ref{thm5} admit an elliptic generalization
and follow by letting $p$ tend to zero
in Theorems~\ref{thm3e} and \ref{thm5e} of the next section.
Before proving the remaining Theorems~\ref{thm2b} and \ref{thm4} we
prepare two simple summation formulae.
\begin{lemma}\label{lem1}
For $c=-abq$ or $c=a^2q/b$ there holds
\begin{equation}\label{phi43sum}
{_4\phi_3}\biggl[\genfrac{}{}{0pt}{}{aq,a^2,b,q^{-n}}
{a,c,a^2bq^{2-n}/c};q,q\biggr]=
\frac{1+aq^n/b}{1+a/b}\frac{(c/a^2q,c/bq;q)_n}{(c,c/a^2bq;q)_n}.
\end{equation}
\end{lemma}
\begin{lemma}\label{lem2}
There holds
\begin{multline}\label{W87}
{_8}W_7(a;b,aq^n/b^{1/2},-aq^n/b^{1/2},q^{-n},-q^{-n};q,q^2) \\
=\frac{(-a/b;q)_{2n}}{(-aq;q)_{2n}}
\frac{(a^2q^2,b;q^2)_n}{(1/b,a^2q^2/b^2;q^2)_n}
\Bigl(\frac{q}{b}\Bigr)^n.
\end{multline}
\end{lemma}

\begin{proof}[Proof of Lemma~\ref{lem1}]
By Sears' transformation for $_4\phi_3$ series \cite[Equation (III.15)]{GR90}
\begin{align*}
\text{LHS}\eqref{phi43sum}
&=\frac{(c/a^2q,c/b;q)_n}{(c,c/a^2bq;q)_n}\,
{_4\phi_3}\biggl[\genfrac{}{}{0pt}{}{a^{-1},b,q^{-1},q^{-n}}
{a,c/a^2q,bq^{1-n}/c};q,q\biggr] \\
&=\frac{(c/a^2q,c/b;q)_n}{(c,c/a^2bq;q)_n}
\Biggl\{1-\frac{a(1-b)(1-q^n)}{b(1-a^2q/c)(1-cq^{n-1}/b)}\biggr\}.
\end{align*}
In general the term within the curly braces does not factor,
but since
\begin{equation*}
\biggl\{\dots\biggr\}=\frac{(1+aq^n/b)(1-c/bq)}{(1+a/b)(1-cq^{n-1}/b)}+
\frac{(1-q^n)(1+c/abq)(1-bc/a^2q)}{(1+b/a)(1-c/a^2q)(1-cq^{n-1}/b)}
\end{equation*}
it certainly does for $c=-abq$ and $c=a^2q/b$, leading to the
right-hand side of \eqref{phi43sum}.
\end{proof}

\begin{proof}[Proof of Lemma~\ref{lem2}]
By Watson's transformation \cite[Equation (III.17)]{GR90}
\begin{equation*}
\text{LHS}\eqref{W87}=
\frac{(aq,b^{1/2};q)_n b^{-n}}{(aq/b,b^{-1/2};q)_n}\,
{_4\phi_3}\biggl[\genfrac{}{}{0pt}{}{b^{1/2}q,b,aq^n/b^{1/2},q^{-n}}
{b^{1/2},-b^{1/2}q^{1-n},-aq^{n+1}};q,q\biggr].
\end{equation*}
The $_4\phi_3$ series on the right can be summed by the $c=-abq$
instance of \eqref{phi43sum}, leading to the desired right-hand side.
\end{proof}

\begin{proof}[Proof of Theorem~\ref{thm2b}]
We write the claim of the theorem as
\begin{align*}
\alpha'_n(a,k;q)&=L_n(a,k;q)\alpha_n(a,m;q), \\
\beta'_n(a,k;q)&=\sum_{r=0}^n N_{n,r}(a,k;q) \beta_r(a,m;q)
\end{align*}
and use the notation \eqref{bMa}.
Then on the one hand
\begin{align*}
\beta'_n(a,k;q)&=\sum_{r=0}^n M_{n,r}(a,k;q)\alpha'_r(a,m;q) \\
&=\sum_{r=0}^n M_{n,r}(a,k;q)L_r(a,k;q)\alpha_r(a,m;q),
\end{align*}
and on the other hand
\begin{align*}
\beta'_n(a,k;q)&=\sum_{s=0}^n N_{n,s}(a,k;q)\beta_s(a,m;q) \\
&=\sum_{s=0}^n N_{n,s}(a,k;q)
\sum_{r=0}^s M_{s,r}(a,m;q) \alpha_r(a,m;q) \\
&=\sum_{r=0}^n \alpha_r(a,m;q)
\sum_{s=r}^n N_{n,s}(a,k;q)M_{s,r}(a,m;q).
\end{align*}
Hence the proof of the theorem boils down to showing that
\begin{equation}\label{NMML}
\sum_{s=r}^n N_{n,s}(a,k;q)M_{s,r}(a,m;q)=M_{n,r}(a,k;q)L_r(a,k;q),
\end{equation}
where $m=a^2/k$.
Using the explicit expressions for $N_{n,s}$ and $M_{r,s}$, 
shifting the summation index $s$ to $s+r$ and carrying out some 
standard manipulations involving $q$-shifted factorials, we obtain
\begin{equation*}
\text{LHS}\eqref{NMML}=
\frac{(m;q)_{2r}}{(aq;q)_{2r}}N_{n,r}(a,k;q)\,
{_4\phi_3}\biggl[\genfrac{}{}{0pt}{}{-\sigma m^{1/2}q^{r+1},mq^{2r},m/a,
q^{-(n-r)}}{-\sigma m^{1/2}q^r,aq^{2r+1},mq^{r-n+1}/k};q,q\biggr].
\end{equation*}
Recalling that $m=a^2/k$ we can sum the $_4\phi_3$ by the $c=a^2q/b$ case of 
Lemma~\ref{lem1}. Thus,
\begin{equation*}
{_4\phi_3}\biggl[\dots\biggr]=
\frac{1-\sigma k^{1/2}q^n}{1-\sigma k^{1/2}q^r}
\frac{(kq^{2r},k/a;q)_{n-r}}{(aq^{2r+1},k/m;q)_{n-r}}
=\frac{(aq;q)_{2r}}{(m;q)_{2r}}
\frac{M_{n,r}(a,k;q)L_r(a,k;q)}{N_{n,r}(a,k;q)}.
\end{equation*}
Putting the above two equations together yields \eqref{NMML}.
\end{proof}

\begin{proof}[Proof of Theorem~\ref{thm4}]
Writing the claim of the theorem as
\begin{align*}
\alpha'_n(a^2,k;q^2)&=L_n(a,k;q)\alpha_n(a,m;q), \\
\beta'_n(a^2,k;q^2)&=\sum_{r=0}^n N_{n,r}(a,k;q) \beta_r(a,m;q)
\end{align*}
and following the proof of Theorem~\ref{thm2b}, we have to show that
\begin{equation}\label{NMML2}
\sum_{s=r}^n N_{n,s}(a,k;q)M_{s,r}(a,m;q)=M_{n,r}(a^2,k;q^2)L_r(a,k;q).
\end{equation}
Using the explicit forms of $N_{n,s}$ and $M_{r,s}$ and
shifting the summation index, we get
\begin{multline*}
\text{LHS}\eqref{NMML2}
=\frac{(m;q)_{2r}}{(aq;q)_{2r}}N_{n,r}(a,k;q) \\
\times
{_8}W_7(mq^{2r};m/a,k^{1/2}q^{n+r},-k^{1/2}q^{n+r},
q^{-(n-r)},-q^{-(n-r)};q;q^2).
\end{multline*}
Since $m=k/a$ the $_8W_7$ can be summed by Lemma~\ref{lem2}. Hence
\begin{multline*}
{_8}W_7(\dots)=
\frac{(-aq^{2r};q)_{2n-2r}}{(-mq^{2r+1};q)_{2n-2r}}
\frac{(m^2q^{4r+2},m/a;q^2)_{n-r}}{(a/m,a^2q^{4r+2};q^2)_{n-r}}
\Bigl(\frac{aq}{m}\Bigr)^{n-r} \\
=\frac{(aq;q)_{2r}}{(m;q)_{2r}}
\frac{M_{n,r}(a^2,k;q^2)L_r(a,k;q)}{N_{n,r}(m^2,k;q^2)}. \qedhere
\end{multline*}
\end{proof}

\section{The elliptic WP Bailey lemma}\label{secell}
We denote by $\theta(z;p)$ the modified Jacobi theta function
\begin{equation*}
\theta(z;p)=\prod_{j=0}^{\infty}(1-zp^j)(1-p^{j+1}/z),\qquad |p|<1,
\end{equation*}
and define the elliptic $q$-shifted factorial by the product
\begin{equation*}
(a;q,p)_n=\prod_{j=0}^{n-1}\theta(aq^j;p).
\end{equation*}
Note that $(a;q,0)_n=(a;q)_n$.
As usual we employ the condensed notation  
\begin{equation*}
(a_1,\dots,a_k;q,p)_n=(a_1;q,p)_n\dots(a_k;q,p)_n.
\end{equation*}

In analogy with the previous section we define the
very-well-poised elliptic hypergeometric series
\begin{equation*}
{_{r+1}V_r}(a_1;a_6,\dots,a_{r+1};q,p)=
\sum_{k=0}^{\infty}\frac{\theta(a_1 q^{2k};p)}{\theta(a_1;p)}
\frac{(a_1,a_6,\dots,a_{r+1};q,p)_k \, q^k}
{(q,a_1q/a_6,\dots,a_1q/a_{r+1};q,p)_k},
\end{equation*}
where for convergence reasons we require the $_{r+1}V_r$ to terminate.
The rationale behind the above labelling of the $_{r+1}V_r$ series is 
that \cite{Spiridonov02a}
\begin{equation*}
\frac{\theta(a q^{2k};p)}{\theta(a;p)}=
\frac{(a^{1/2}q,-a^{1/2}q,(a/p)^{1/2}q,-(ap)^{1/2}q;q,p)_k}
{(a^{1/2},-a^{1/2},(ap)^{1/2},-(a/p)^{1/2};q,p)_k}(-q)^k.
\end{equation*}
Hence, provided none of the parameters $a_i$ depends on $p$
(see Theorems~\ref{thm1413b} and \ref{thm1413c} for examples of 
such $p$-dependence) we get
\begin{equation}\label{VW}
{_{r+1}V_r}(a_1;a_6,\dots,a_{r+1};q,0)=
{_{r-1}W_{r-2}}(a_1;a_6,\dots,a_{r+1};q,q).
\end{equation}
A $_{r+1}V_r$ series is called balanced if
$a_6\cdots a_{r+1}q=(a_1 q)^{(r-5)/2}$. 
All known identities for elliptic hypergeometric series are
both balanced and very-well-poised.
One such identity that will be applied on a number of occasions is
the elliptic version of Jackson's $_8\phi_7$ sum due to
Frenkel and Turaev \cite[Theorem 5.5.2]{FT97};
\begin{equation}\label{V109}
{_{10}V_9}(a;b,c,d,e,q^{-n};q,p)=
\frac{(aq,aq/bc,aq/bd,aq/cd;q,p)_n}{(aq/b,aq/c,aq/d,aq/bcd;q,p)_n},
\end{equation}
where $bcde=a^2 q^{n+1}$. For a more extensive introduction to 
elliptic hypergeometric series we refer the reader to 
\cite{vDS00,vDS01a,vDS01b,vDS02,FT97,KNR03,Rosengren01,Rosengren02,RS03,
Spiridonov02a,Spiridonov02,Spiridonov03,SZ99,W02}.

\vspace{4mm}

As our first application of elliptic hypergeometric series
we state the following analogue of the WP Bailey transform.
\begin{lemma}[Elliptic WP Bailey transform]
For $a$ and $k$ indeterminates the following two equations are equivalent:
\begin{subequations}\label{WPBTe}
\begin{align}\label{abe}
\beta_n(a,k;q,p)&=\sum_{r=0}^n \frac{(k/a;q,p)_{n-r}}{(q;q,p)_{n-r}}
\frac{(k;q,p)_{n+r}}{(aq;q,p)_{n+r}}\,\alpha_r(a,k;q,p), \\
\alpha_n(a,k;q,p)&=
\frac{\theta(aq^{2n};p)}{\theta(a;p)} 
\sum_{r=0}^n\frac{\theta(kq^{2r};p)}{\theta(k;p)}
\frac{(a/k;q,p)_{n-r}}{(q;q,p)_{n-r}}
\label{bae} \\
& \qquad \qquad \qquad \qquad \times \frac{(a;q,p)_{n+r}}{(kq;q,p)_{n+r}}
\Bigl(\frac{k}{a}\Bigr)^{n-r}\beta_r(a,k;q,p). \notag
\end{align}
\end{subequations}
\end{lemma}
\begin{proof}
The proof is an almost exact copy of the proof of Lemma~\ref{lemT}.
Writing the two transformations of the lemma as
\begin{subequations}
\begin{align}\label{bMae}
\beta_n(a,k;q,p)&=\sum_{r=0}^n M_{n,r}(a,k;q,p) \alpha_r(a,k;q,p), \\
\alpha_n(a,k;q,p)&=\sum_{r=0}^n \tilde{M}_{n,r}(a,k;q,p) \beta_r(a,k;q,p),
\end{align}
\end{subequations}
we have to show the inverse relation
\begin{equation*}
\sum_{s=r}^n \tilde{M}_{n,s}(a,k;q,p) M_{s,r}(a,k;q,p)=\delta_{n,r}.
\end{equation*}
This is equivalent to showing that
\begin{equation*}
{_8}V_7(kq^{2s};k/a,aq^{n+r},q^{-(n-r)};q,p)=\delta_{n,r},
\end{equation*}
and readily follows from the elliptic Jackson sum \eqref{V109}.
\end{proof}

The next result is Spiridonov's elliptic version of 
Andrews' Theorem~\ref{thm1} \cite[Theorem 4.3]{Spiridonov02}.
\begin{theorem}\label{thm1e}
If $(\alpha(a,k;q,p),\beta(a,k;q,p))$ is an elliptic WP Bailey pair, 
then so is the pair $(\alpha'(a,k;q,p),\beta'(a,k;q,p)$ given by
\begin{align*}
\alpha'_n(a,k;q,p)&=\frac{(b,c;q,p)_n}{(aq/b,aq/c;q,p)_n}
\Bigl(\frac{k}{m}\Bigr)^n \alpha_n(a,m;q,p), \\
\beta'_n(a,k;q,p)&=\frac{(mq/b,mq/c;q,p)_n}
{(aq/b,aq/c;q,p)_n}\sum_{r=0}^n \frac{\theta(mq^{2r};p)}{\theta(m;p)}
\frac{(b,c;q,p)_r}{(mq/b,mq/c;q,p)_r} \\
& \qquad \qquad \qquad \qquad \times
\frac{(k/m;q,p)_{n-r}}{(q;q,p)_{n-r}}
\frac{(k;q,p)_{n+r}}{(mq;q,p)_{n+r}} 
\Bigl(\frac{k}{m}\Bigr)^r \beta_r(a,m;q,p),
\end{align*}
where $m=bck/aq$.
\end{theorem}
To this we can add two more transformations for elliptic WP Bailey pairs.
The first generalizes Theorem~\ref{thm3}.
\begin{theorem}\label{thm3e}
If $(\alpha(a,k;q,p),\beta(a,k;q,p))$ is an elliptic WP Bailey pair, 
then so is the pair $(\alpha'(a,k;q,p),\beta'(a,k;q,p)$ given by
\begin{align*}
\alpha'_n(a^2,k;q^2,p^2)&=\alpha_n(a,m;q,p), \\
\beta'_n(a^2,k;q^2,p^2)&=
\frac{(-mq;q,p)_{2n}}{(-aq;q,p)_{2n}}
\sum_{r=0}^n \frac{\theta(mq^{2r};p)}{\theta(m;p)}
\frac{(k/m^2;q^2,p^2)_{n-r}}{(q^2;q^2,p^2)_{n-r}} \\
&\qquad \qquad \qquad \qquad \quad \times
\frac{(k;q^2,p^2)_{n+r}}{(m^2q^2;q^2,p^2)_{n+r}}
\Bigl(\frac{m}{a}\Bigr)^{n-r} \beta_r(a,m;q,p),
\end{align*}
where $m=k/aq$.
\end{theorem}
The second result provides an elliptic extension of Theorem~\ref{thm5}.
\begin{theorem}\label{thm5e}
If $(\alpha(a,k;q,p),\beta(a,k;q,p))$ is an elliptic WP Bailey pair,
then so is the pair $(\alpha'(a,k;q,p),\beta'(a,k;q,p)$ given by
\begin{align*}
\alpha'_{2n}(a,k;q,p)&=\alpha_n(a,m;q^2,p), \qquad  
\alpha'_{2n+1}(a,k;q,p)=0, \\
\beta'_n(a,k;q,p)&=\frac{(mq;q^2,p)_n}{(aq;q^2,p)_n}
\sum_{r=0}^{\lfloor n/2\rfloor}
\frac{\theta(mq^{4r};p)}{\theta(m;p)}
\frac{(k/m;q,p)_{n-2r}}{(q;q,p)_{n-2r}} \\
& \qquad\qquad\qquad\qquad\quad \times
\frac{(k;q,p)_{n+2r}}{(mq;q,p)_{n+2r}}
\Bigl(-\frac{k}{a}\Bigr)^{n-2r} 
\beta_r(a,m;q^2,p),
\end{align*}
where $m=k^2/a$.
\end{theorem}

\begin{proof}[Proof of Theorem~\ref{thm3e}]
Writing the claim of the theorem as
\begin{align*}
\alpha'_n(a^2,k;q^2,p^2)&=\alpha_n(a,m;q,p), \\
\beta'_n(a^2,k;q^2,p^2)&=\sum_{r=0}^n N_{n,r}(a,k;q,p) \beta_r(a,m;q,p)
\end{align*}
and using the notation of equation \eqref{bMae}, we need to show that
\begin{equation}\label{NMMe}
\sum_{s=r}^n N_{n,s}(a,k;q,p)M_{s,r}(a,m;q,p)=M_{n,r}(a^2,k;q^2,p^2),
\end{equation}
where $m=k/aq$.
From the explicit expressions for $N_{n,s}$ and $M_{r,s}$, 
and simple elliptic $q$-factorial relations such as
\begin{subequations}\label{sm}
\begin{align}
(a;q,p)_{n+k}&=(a;q,p)_n (aq^n;q,p)_k \\
(a;q,p)_{n-k}&=(a;q,p)_n(-q^{1-n}/a)^k q^{\binom{k}{2}}/(q^{1-n}/a;q,p)_k \\
(a^2;q^2,p^2)_n&=(a,-a;q,p)_n,
\end{align}
\end{subequations}
we obtain
\begin{multline*}
\text{LHS}\eqref{NMMe}
=\frac{(m;q,p)_{2r}}{(aq;q,p)_{2r}}N_{n,r}(m^2,k;q^2,p^2) \\
\times
{_{10}V_9}(mq^{2r};m/a,k^{1/2}q^{n+r},-k^{1/2}q^{n+r},
q^{-(n-r)},-q^{-(n-r)};q,p).
\end{multline*}
Since $m=k/aq$ the $_{10}V_9$ can be summed by \eqref{V109} to yield
\begin{multline*}
{_{10}V_9}(\dots)=
\frac{(kq^{2r}/a,aq^{r-n+1}/k^{1/2},-aq^{r-n+1}/k^{1/2},
-q^{-2n}/a;q,p)_{n-r}}
{(aq^{2r+1},k^{1/2}q^{r-n}/a,-k^{1/2}q^{r-n}/a,-aq^{-2n+1}/k;q,p)_{n-r}} \\
=\frac{(aq;q,p)_{2r}}{(m;q,p)_{2r}}
\frac{M_{n,r}(a^2,k;q^2,p^2)}{N_{n,r}(m^2,k;q^2,p^2)}. \qedhere
\end{multline*}
\end{proof}

\begin{proof}[Proof of Theorem~\ref{thm5e}]
We write the claim of the theorem as
\begin{align*}
\alpha'_{2n}(a,k;q,p)&=\alpha_n(a,m;q^2,p), \qquad 
\alpha'_{2n+1}(a,k;q,p)=0, \\
\beta'_n(a,k;q,p)&=\sum_{r=0}^{\lfloor n/2\rfloor}
N_{n,r}(a,k;q,p) \beta_r(a,m;q^2,p)
\end{align*}
and again use the notation \eqref{bMae}.
We then need to show that
\begin{equation}\label{NMM2e}
\sum_{s=r}^{\lfloor n/2\rfloor} N_{n,s}(a,k;q,p)
M_{s,r}(a,m;q^2,p)=M_{n,2r}(a,k;q,p),
\end{equation}
where $m=k^2/a$.
From \eqref{sm} and $(a;q,p)_{2n}=(a,aq;q^2,p)_n$ it follows that
\begin{multline*}
\text{LHS}\eqref{NMM2e}
=\frac{(m;q^2,p)_{2r}}{(aq^2;q^2,p)_{2r}} N_{n,r}(m,k;q,p) \\
\times
{_{10}V_9}(mq^{4r};m/a,kq^{n+2r},kq^{n+2r+1},
q^{-(n-2r)},q^{1-(n-2r)};q^2,p).
\end{multline*}
Since $m=k^2/a$ the ${_{10}V_9}$ can be summed by \eqref{V109}
upon distinguishing between $n$ even and $n$ odd.
After some simplifications one finds that irrespective of
the parity of $n$
\begin{equation*}
{_{10}V_9}(\dots)=\frac{(aq^2;q^2,p)_{2r}}{(m;q^2,p)_{2r}}
\frac{M_{n,2r}(a,k;q,p)}{N_{n,r}(m,k;q,p)},
\end{equation*}
which completes the proof.
\end{proof}

\section{Applications}\label{secapl}
As mentioned in the introduction, all of the WP and elliptic WP
Bailey pairs found in \cite{Andrews01,AB02,Spiridonov02}
readily follow by application of the Theorems of the previous
two sections.
To illustrate this we take $\beta(a,k;q)=\delta_{n,0}$ in the WP Bailey
transform \eqref{WPBT}, yielding the unit WP Bailey pair
\begin{subequations}\label{UBP}
\begin{align}
\alpha_n(a,k;q)&=\frac{1-aq^{2n}}{1-a}
\frac{(a,a/k;q)_n}{(q,kq;q)_n}\Bigl(\frac{k}{a}\Bigr)^n, \\
\beta_n(a,k;q)&=\delta_{n,0}.
\end{align}
\end{subequations}
Then applying Theorem~\ref{thm2b} results in the (corrected)
WP Bailey pair given by \cite[Equations (3.5) and (3.6)]{AB02},
applying Theorem~\ref{thm5} results in the (corrected)
WP Bailey pair given by \cite[Equations (3.11) and (3.12)]{AB02},
applying Theorem~\ref{thm3} results in the 
WP Bailey pair given by \cite[Equations (4.5) and (4.6)]{AB02},
and, finally,
applying Theorem~\ref{thm4} results in the (corrected)  
WP Bailey pair given by \cite[Equations (4.7) and (4.8)]{AB02}.
These last two cited WP Bailey pairs were first found by
Bressoud in \cite{Bressoud81}.

In much the same way, from the elliptic WP Bailey transform \eqref{WPBTe} one
immediately infers the elliptic unit WP Bailey pair
\begin{subequations}\label{UBPe}
\begin{align}
\alpha_n(a,k;q,p)&=\frac{\theta(aq^{2n};p)}{\theta(a;p)}
\frac{(a,a/k;q,p)_n}{(q,kq;q,p)_n}\Bigl(\frac{k}{a}\Bigr)^n, \\
\beta_n(a,k;q,p)&=\delta_{n,0}.
\end{align}
\end{subequations}
When inserted in Theorem~\ref{thm3e} this yields
\begin{align*}
\alpha_n(a^2,k;q^2,p^2)&=
\frac{\theta(aq^{2n};p)}{\theta(a;p)}
\frac{(a,a^2q/k;q,p)_n}{(q,k/a;q,p)_n}\Bigl(\frac{k}{a^2q}\Bigr)^n, \\
\beta_n(a^2,k;q^2,p^2)&=
\frac{(-k/a;q,p)_{2n}}{(-aq;q,p)_{2n}}
\frac{(k,a^2q^2/k;q^2,p^2)_n}{(q^2,k^2/a^2;q^2,p^2)_n}
\Bigl(\frac{m}{a}\Bigr)^n, 
\end{align*}
which is equivalent to 
\cite[Equations (5.2) and (5.3)]{Spiridonov02} of Spiridonov.
By further iterating this pair using Theorem~\ref{thm1e} Spiridonov obtained
the following transformation formula for elliptic hypergeometric series
\cite[Theorem 5.1]{Spiridonov02}.
\begin{theorem}\label{thm1413b}
For $m=bck/a^2q^2$ and $d=-m/a$ there holds
\begin{align*}
{_{14}V_{13}}&(a;a^2q/m,b^{1/2},-b^{1/2},c^{1/2},-c^{1/2},
k^{1/2}q^n,-k^{1/2}q^n,q^{-n},-q^{-n};q,p) \\
&=\frac{(a^2q^2,k/m,mq^2/b,mq^2/c;q^2,p^2)_n}
{(mq^2,k/a^2,a^2q^2/b,a^2q^2/c;q^2,p^2)_n} \\
& \quad \qquad \times {_{14}V_{13}}(m;a^2q^2/m,d,dq,d/p,dqp,b,c,
kq^{2n},q^{-2n};q^2,p^2).
\end{align*}
\end{theorem}
Interestingly, some of the parameters in the $_{14}V_{13}$ series on
the right depend on $p$. Therefore \eqref{VW} does not apply,
and in the $p\to 0$ limit the above identity ceases to be balanced.
Indeed in this limit one finds
\begin{align*}
{_{12}W_{11}}&(a;a^2q/m,b^{1/2},-b^{1/2},c^{1/2},-c^{1/2},
k^{1/2}q^n,-k^{1/2}q^n,q^{-n},-q^{-n};q,q) \\
&\quad =\frac{(a^2q^2,k/m,mq^2/b,mq^2/c;q^2)_n}
{(mq^2,k/a^2,a^2q^2/b,a^2q^2/c;q^2)_n} \\
& \qquad \qquad \qquad \times
{_{10}W_9}(m;a^2q^2/m,d,dq,b,c,kq^{2n},q^{-2n};q^2,mq/a^2)
\end{align*}
of Andrews and Berkovich \cite[Equation (4.9)]{AB02}.
What we will now show is that thanks to Theorem~\ref{thm5e} the 
above theorem has the following companion.
\begin{theorem}\label{thm1413c}
For $m=bck/aq$ and $d=\pm m(q/a)^{1/2}$ there holds
\begin{align*}
{_{14}V_{13}}&(a;a^2/m^2,b,bq,c,cq,kq^n,kq^{n+1},q^{-n},q^{1-n};q^2,p) \\
&=\frac{(aq,k/m,mq/b,mq/c;q,p)_n}{(mq,k/a,aq/b,aq/c;q,p)_n}  \\
&\qquad\quad \times
{_{14}V_{13}}(m;a/m,d,-d,
dp^{1/2},-d/p^{1/2},b,c,kq^n,q^{-n};q,p).
\end{align*}
\end{theorem}
This provides a second example of an identity that fails to be balanced 
after taking $p$ to zero. In this limit
\begin{multline*}
{_{12}W_{11}}(a;a^2/m^2,b,bq,c,cq,kq^n,kq^{n+1},q^{-n},q^{1-n};q^2,q^2) \\
=\frac{(aq,k/m,mq/b,mq/c;q)_n}{(mq,k/a,aq/b,aq/c;q)_n} \,
{_{10}W_9}(m;a/m,d,-d,b,c,kq^n,q^{-n};q,-mq/a)
\end{multline*}
again due to Andrews and Berkovich \cite[Equation (3.13)]{AB02}.

\begin{proof}[Proof of Theorem~\ref{thm1413c}]
Substituting the elliptic unit WP Bailey pair \eqref{UBPe}
in Theorem~\ref{thm5e} results in the new pair
\begin{align*}
\alpha_{2n}(a,k;q,p)&=
\frac{\theta(aq^{4n};p)}{\theta(a;p)}
\frac{(a,a^2/k^2;q^2,p)_n}{(q^2,k^2q^2/a;q^2,p)_n}
\Bigl(\frac{k}{a}\Bigr)^{2n}, \\
\alpha_{2n+1}(a,k;q,p)&=0, \\
\beta_n(a,k;q,p)&=\frac{(k^2q/a;q^2,p)_n}{(aq;q^2,p)_n}
\frac{(k,a/k;q,p)_n}{(q,k^2q/a;q,p)_n}\Bigl(-\frac{k}{a}\Bigr)^n.
\end{align*}
Next applying Theorem~\ref{thm1e} and manipulating some of the
elliptic $q$-shifted factorials using \eqref{sm} yields
\begin{align*}
\alpha_{2n}(a,k;q,p)&=
\frac{\theta(aq^{4n};p)}{\theta(a;p)}
\frac{(a,a^2/m^2;q^2,p)_n}{(q^2,m^2q^2/a;q^2,p)_n}
\frac{(b,c;q,p)_{2n}}{(aq/b,aq/c;q,p)_{2n}}
\Bigl(\frac{k}{a}\Bigr)^{2n},  \\
\alpha_{2n+1}(a,k;q,p)&=0, \\
\beta_n(a,k;q,p)&=
\frac{(k,k/m,bk/a,ck/a;q,p)_n}{(q,mq,aq/b,aq/c;q,p)_n} \\
&\quad \times \sum_{r=0}^n \frac{\theta(mq^{2r};p)}{\theta(m;p)}
\frac{(m^2q/a;q^2,p)_r}{(aq;q^2,p)_r}
\Bigl(-\frac{mq}{a}\Bigr)^r  \\
& \qquad\qquad \times \frac{(m,b,c,a/m,kq^n,q^{-n};q,p)_r}
{(q,mq/b,mq/c,m^2q/a,mq^{1-n}/k,mq^{n+1};q,p)_r},
\end{align*}
where $m=bck/aq$.
By taking the ratio of 
\begin{equation*}
\theta(a^2;p)
=\theta(a;p)\theta(-a;p)\theta(ap^{1/2};p)\theta(-a/p^{1/2};p)
p^{1/2}/a
\end{equation*}
and the equation obtained by replacing $a$ by $-b$ it readily follows that
\begin{equation*}
\frac{(a^2;q^2,p)_n}{(b^2;q^2,p)_n}
=\frac{(a,-a,ap^{1/2},-a/p^{1/2};q,p)_n}
{(b,-b,b/p^{1/2},-bp^{1/2};q,p)_n}\Bigl(-\frac{b}{a}\Bigr)^n.
\end{equation*}
Defining $d$ by $d^2=m^2q/a$ one thus finds
\begin{equation*}
\frac{(m^2q/a;q^2,p)_r}{(aq;q^2,p)_r}
\Bigl(-\frac{m}{a}\Bigr)^r
=\frac{(d,-d,dp^{1/2},-d/p^{1/2};q,p)_r}
{(mq/d,-mq/d,mq/dp^{1/2},-mqp^{1/2}/d;q^2,p)_r}.
\end{equation*}
Using this in the expression for $\beta_n$ and $(a;q,p)_{2n}=
(a,aq;q^2,p)_n$ in the expression for $\alpha_n$, and
substituting the above Bailey pair in \eqref{abe}
yields the $_{14}V_{13}$ transformation of the theorem.
\end{proof}

Many more transformations for basic and elliptic hypergeometric series
can be derived along the lines of the above proof. 
Most of the simpler results are either known or variations of known
identities.
To give just two more examples, iterating the unit Bailey pair 
\eqref{UBP} using first
Theorem~\ref{thm1} and then Theorem~\ref{thm2b} or Theorem~\ref{thm4}
yields
\begin{multline*}
{_{12}W_{11}}(a;b,c,kq/bc,(mq)^{1/2},m^{1/2}q,-m^{1/2},-(mq)^{1/2},
kq^n,q^{-n};q,q) \\
=\frac{1+k^{1/2}}{1+k^{1/2}q^n}
\frac{(aq,k/m;q)_n}{(k,k/a;q)_n}\,
{_6\phi_5}\biggl[\genfrac{}{}{0pt}{}
{m,m^{1/2}q,bm/a,cm/a,aq/bc,q^{-n}}
{m^{1/2},aq/b,aq/c,bcm/a,mq^{1-n}/k};q,q\biggr],
\end{multline*}
with $m=a^2/k$, and
\begin{align*}
{_{12}W_{11}}&(a;b,c,a^2q/bcm,(-a)^{1/2}q,-(-a)^{1/2}q,
k^{1/2}q^n,-k^{1/2}q^n,q^{-n},-q^{-n};q,q) \\
&=\frac{(-mq;q)_{2n}}{(-a;q)_{2n}}
\frac{(a^2q^2,k/m^2;q^2)_n}{(m^2q^2,k/a^2;q^2)_n} 
\Bigl(\frac{m}{aq}\Bigr)^n \\
& \qquad \times
{_{10}W_9}(m;bm/a,cm/a,aq/bc,k^{1/2}q^n,-k^{1/2}q^n,q^{-n},-q^{-n};q,q^2),
\end{align*}
with $m=k/a$.
The first of these formulas is similar to a nearly-poised
transformation of Bailey \cite[Equation (III.25)]{GR90} and the second
formula is similar to a special case of Bailey's $_{10}\phi_9$ transformation
\cite[Equation (III.28)]{GR90}.

\section{Note added in proof}
After submission of this paper we discovered many more transformations
for basic and elliptic WP Bailey pairs.
The most important two, which hold at the elliptic level,
are given below.
\begin{theorem}\label{thmnew1}
If $(\alpha(a,k;q,p),\beta(a,k;q,p))$ is an elliptic WP Bailey pair,
then so is the pair $(\alpha'(a,k;q,p),\beta'(a,k;q,p)$ given by
\begin{align*}
\alpha'_n(a,k;q,p)&=\frac{(mq;q^2,p)_n}{(kq;q^2,p)_n}
\Bigl(-\frac{a}{m}\Bigr)^n \alpha_n(a,m;q,p), \\
\beta'_n(a,k;q,p)&=
\sum_{\substack{r=0 \\[0.3mm] r\equiv n\:(2)}}^n 
\frac{\theta(mq^{2r};p)}{\theta(m;p)} 
\frac{(k/m;q^2,p)_{(n-r)/2}}{(q^2;q^2,p)_{(n-r)/2}} \\[-1.5mm]
& \qquad \qquad \qquad \qquad \times
\frac{(k;q^2,p)_{(n+r)/2}}{(mq^2;q^2,p)_{(n+r)/2}} 
\Bigl(-\frac{a}{m}\Bigr)^r \beta_r(a,m;q,p),
\end{align*}
where $m=a^2/k$.
\end{theorem}
\begin{theorem}\label{thmnew2}
If $(\alpha(a,k;q,p),\beta(a,k;q,p))$ is an elliptic WP Bailey pair,
then so is the pair $(\alpha'(a,k;q,p),\beta'(a,k;q,p)$ given by
\begin{align*}
\alpha'_n(a,k^2;q^2,p^2)&=\frac{(-mq;q,p)_{2n}}{(-kq;q,p)_{2n}}
\Bigl(\frac{a}{m^2q}\Bigr)^n \alpha_n(a,m^2;q^2,p^2), \\
\beta'_n(a,k^2;q^2,p^2)&=
\frac{\theta(-k;p)}{\theta(-kq^{2n};p)}\,q^{-n} 
\sum_{r=0}^n 
\frac{\theta(m^2q^{4r};p^2)}{\theta(m^2;p^2)} 
\frac{(k/m;q,p)_{n-r}}{(q;q,p)_{n-r}} \\
& \qquad \qquad \qquad \qquad \quad \times
\frac{(k;q,p)_{n+r}}{(mq;q,p)_{n+r}} 
\Bigl(\frac{a}{m^2}\Bigr)^r \beta_r(a,m^2;q^2,p^2),
\end{align*}
where $m=a/kq$.
\end{theorem}
Especially this last theorem is remarkable because, unlike all the other
transformations for elliptic WP Bailey pairs, its proof does not
rely on the elliptic Jackson sum \eqref{V109} 
but on the new elliptic summation
\begin{multline*}
\sum_{k=0}^n \frac{\theta(a^2q^{4k};p^2)}{\theta(a^2;p^2)}
\frac{(a^2,b;q^2,p^2)_k}{(q^2,a^2 q^2/b;q^2,p^2)_k}
\frac{(a q^{n-1}/b,q^{-n};q,p)_k}{(bq^{2-n},aq^{n+1};q,p)_k}\,q^{2k} \\
=\frac{\theta(-aq^{2n-1}/b;p)}{\theta(-a/bq;p)}
\frac{(aq,-a/bq;q,p)_n}{(-q,1/bq;q,p)_n}
\frac{(1/bq^2;q^2,p^2)_n}{(a^2q^2/b;q^2,p^2)_n}\, q^n.
\end{multline*}
When $p$ tends to zero this reduces to a bibasic summation 
of Nassrallah and Rahman \cite[Corollary 4]{NR81} (see also
\cite[Equation (3.10.8)]{GR90}).

An appealing example of a new result that follows from 
Theorem~\ref{thmnew2} (combined with \eqref{UBPe} and Theorem~\ref{thm1e}) 
is
\begin{align*}
&\sum_{k=0}^n \frac{\theta(a^2q^{4k};p^2)}{\theta(a^2;p^2)}
\frac{(a^2,b,c,d;q^2,p^2)_k}{(q^2,a^2q^2/b,a^2q^2/c,a^2q^2/d;q^2,p^2)_k}
\frac{(e q^n,q^{-n};q,p)_k}
{(aq^{1-n}/e,aq^{n+1};q,p)_k}\,q^{2k} \\
& =\frac{\theta(-e q^{2n};p)}{\theta(-e;p)}
\frac{(-e,aq;q,p)_n}{(-q,e/a;q,p)_n}
\frac{(e/aq;q^2,p^2)_n}{(\lambda q^2;q^2,p^2)_n}\, q^n \\
&\quad \times 
{_{14}V_{13}}(\lambda,-aq,-aq^2,-aq/p,-aq^2p,
\lambda b/a^2,\lambda c/a^2,\lambda d/a^2,e^2 q^{2n},q^{-2n};q^2,p^2),
\end{align*}
where $\lambda=a^4 q^2/bcd$ and $e=\lambda/a q$.
Again the $p$-dependence of some of the parameters of the 
$_{14}V_{13}$ on the right is to be noted.
In the limit when $p$ tends to zero the above transformation
simplifies to \cite[Equation (4.24)]{NR81} (see also
\cite[Exercise 3.13]{GR90}.

\vspace{4mm}

As a final comment we wish to remark that not all of the
(elliptic) WP Bailey pairs one can find seem
to follow from the (elliptic) unit WP Bailey pair by
transformations of the type discussed in this paper.
One example is
\begin{align*}
\alpha_n(a,k;q^2,p)&=
\frac{\theta(aq^{4n};p)}{\theta(a;p)}
\frac{(a,m^2q/k,b,a/b;q^2,p)_n}{(q^2,akq/m^2,aq^2/b,bq^2;q^2,p)_n}
\frac{(aq/m;q,p)_{2n}}{(m;q,p)_{2n}}
\Bigl(\frac{k}{a}\Bigr)^n, \\
\beta_n(a,k;q^2,p)&=\frac{(m^2q/a;q^2,p)_n}{(akq/m^2;q^2,p)_n}
\sum_{r=0}^n \frac{\theta(mq^{3r};p)}{\theta(m;p)}
\frac{(aq/m,bm/a,m/b;q,p)_r}{(m^2q/a,aq^2/b,bq^2;q^2,p)_r} \\
&\qquad \times
\frac{(m^2q/k;q^2,p)_r}{(k/m;q,p)_r}
\frac{(k/m;q,p)_{2n-r}}{(q^2;q^2,p)_{n-r}}
\frac{(k;q^2,p)_{n+r}}{(mq;q,p)_{2n+r}} \,
q^{\binom{r}{2}}\Bigl(\frac{k}{m}\Bigr)^r, 
\end{align*}
which implies the transformation \cite[Theorem 4.2]{W02}
\begin{align*}
\sum_{r=0}^n & \frac{\theta(mq^{3r};p)}{\theta(m;p)}
\frac{(aq/m,bm/a,m/b;q,p)_r}{(m^2q/a,aq^2/b,bq^2;q^2,p)_r}
\frac{(m^2q/k,kq^{2n},q^{-2n};q^2,p)_r}
{(k/m,mq^{1-2n}/k,mq^{2n+1};q,p)_r} \, q^r \\
&\qquad =\frac{(k/a,akq/m^2;q^2,p)_n}{(aq^2,m^2q/a;q^2,p)_n}
\frac{(mq;q,p)_{2n}}{(k/m;q,p)_{2n}} \\
& \qquad \qquad \quad \times
{_{12}V_{11}}(a;b,a/b,m^2q/k,aq/m,aq^2/m,kq^{2n},q^{-2n};q^2,p).
\end{align*}
In order to establish the above elliptic WP Bailey pair
(and similar pairs not given here) one has to extend some of our 
earlier results to bibasic series.
A first step is to define a bibasic elliptic WP Bailey pair by
\begin{equation*}
B^{(i)}_n(a,k;q,p)=\sum_{r=0}^n 
\frac{(k/a;q,p)_{n-ir}}{(q^i;q^i,p)_{n-r}}
\frac{(k;q,p)_{n+ir}}{(aq^i;q^i,p)_{n+r}}\,A^{(i)}_r(a,k;q,p),
\end{equation*}
with $i$ a positive integer.
Note that $(A^{(1)}_n,B^{(1)}_n)=(\alpha_n,\beta_n)$.
The next and non-trivial step is to relate a pair
$(A^{(i)}_n,B^{(i)}_n)$ with $i\geq 2$ to $(\alpha_n,\beta_n)$.
Two examples of such relations can be stated as follows.
\begin{theorem}\label{thm2to1}
There holds
\begin{align*}
\alpha_n(a,k;q^2,p)&=\frac{(m^2q/k;q^2,p)_n}{(akq/m^2;q^2,p)_n}
q^{n^2}\Bigl(-\frac{ak}{m^2}\Bigr)^n A^{(2)}_n(a,m;q,p), \\
\beta_n(a,k;q^2,p)&=\frac{(m^2q/a;q^2,p)_n}{(akq/m^2;q^2,p)_n}
\sum_{r=0}^n \frac{\theta(mq^{3r};p)}{\theta(m;p)}
\frac{(aq/m;q,p)_r}{(m^2q/a;q^2,p)_r}
\frac{(m^2q/k;q^2,p)_r}{(k/m;q,p)_r} \\
& \qquad \times 
\frac{(k/m;q,p)_{2n-r}}{(q^2;q^2,p)_{n-r}}
\frac{(k;q^2,p)_{n+r}}{(mq;q,p)_{2n+r}} \,
q^{\binom{r}{2}}\Bigl(\frac{k}{m}\Bigr)^r B^{(2)}_r(a,m;q,p).
\end{align*}
\end{theorem}
\begin{theorem}
For $m^3=ak$ there holds
\begin{align*}
\alpha_n(a,k;q^3,p)&=q^{3n^2}a^r A^{(3)}_n(a,m;q,p), \\
\beta_n(a,k;q^3,p)&=
\sum_{r=0}^n \frac{\theta(mq^{4r};p)}{\theta(m;p)}
\frac{(aq/m;q,p)_{2r}}{(m^2/a;q,p)_{2r}} \\
& \qquad \times 
\frac{(k/m;q,p)_{3n-r}}{(q^3;q^3,p)_{n-r}}
\frac{(k;q^3,p)_{n+r}}{(mq;q,p)_{3n+r}} \,
q^{2\binom{r}{2}}\Bigl(\frac{k}{m}\Bigr)^r B^{(3)}_r(a,m;q,p).
\end{align*}
\end{theorem}
The proofs hinge on \cite[Corollaries 4.4 and 4.5]{W02}.
As third and final step in the proof of the elliptic WP Bailey pair
of our example
we note that \cite[Theorem 4.1]{W02} is equivalent to
\begin{align*}
A^{(i)}_n(a,k;q,p)&=\frac{\theta(aq^{2in};p)}{\theta(a;p)}
\frac{(a,b,a/b;q^i,p)_n}{(q^i,aq^i/b,bq^i;q^i,p)_n} \\
& \qquad \qquad \qquad \times
\frac{(aq/k;q,p)_{in}}{(k;q,p)_{in}}\,
(-1)^n q^{-\binom{i}{2}n^2} \Bigl(-\frac{k}{a}\Bigr)^{in},\\
B^{(i)}_n(a,k;q,p)&=\frac{(bk/a,k/b;q,p)_n}{(aq^i/b,bq^i;q^i,p)_n}.
\end{align*}
Substituting this in Theorem~\ref{thm2to1} concludes our derivation.

We leave it as an open problem to find transformations between
$(A^{(i)}_n,B^{(i)}_n)$ and $(A^{(j)}_n,B^{(j)}_n)$ for arbitrary 
$i$ and $j$.

\bibliographystyle{amsplain}

\end{document}